\documentclass[twoside,11pt,reqno]{amsart}

\usepackage{amsmath,amsthm,amssymb,amstext,amsfonts,amscd}
\usepackage{graphicx}
\usepackage{multirow}
\usepackage{url}

\numberwithin{equation}{section}

\begin{document}
\title[{\bf Study of dynamical symmetrietry algebra of $\Psi_{2}$--Humbert function}]
{Study of dynamical symmetrietry algebra of $\Psi_{2}$--Humbert function}

\author[{\bf A. Shehata and D. Kumar}]{\bf Ayman Shehata and Dinesh Kumar$^{*}$}

\address{A. Shehata: Department of Mathematics, Faculty of Science, Assiut University, Assiut--71516,  EGYPT}
\email{drshehata2006@yahoo.com; aymanshehata@science.aun.edu.eg}

\address{D. Kumar: Department of Applied Sciences, College of Agriculture--Jodhpur,  Agriculture University Jodhpur, Jodhpur--342304, INDIA}
\email{dinesh\_dino03@yahoo.com}

\bigskip

\thanks{$^*$ Corresponding author}



\begin{abstract}
The study is devoted to the construction of dynamical symmetry algebra of confluent hypergeometric function $\;_{1}F_{1}$ and $\Psi_{2}$--Humbert function and to derive some generating relations and reduction formulas for  $\;_{1}F_{1}$ and $\Psi_{2}$ functions.
\end{abstract}
\maketitle

\section {Introduction and Preliminaries}\label{sec 1}
Special functions in mathematical physics, often viewed as solutions to partial differential equations such as the Laplace and diffusion equations, have been explored by various authors through diverse approaches. Beyond their applications in physics and engineering, the theoretical aspects of these functions have captivated mathematicians for over a century. The author’s ongoing investigation into the Lie algebraic structure of hypergeometric-type functions is further detailed in \cite{wm1, wm3}.

Dynamical symmetry algebra of hypergeometric function was constructed by Milier \cite{wm1, wm2, wm8, wm5} and its further use was presented by Agarwal and Jain \cite{ar} to find generating functions for Jcobi polynomials. Miller \cite{wm6, wm8, wm9} discussed the Lie theory of the Appell function $F_{1}$, Lauricella functions $F_{D}$ and Meijers’s $G$-Function, which has proven to be its dynamical symmetry algebra. In \cite{lw} Louis Weiner investigated the group-theoretic origin of certain generating functions. In \cite{ca} Carlitz introduced some reduction formulas for generalized hypergeometric functions. Kishan Sharma and Renu Jain \cite{sj3, sj4, sj5, sj6}  derived  the dynamical symmetry algebra of  basic confluent hypergeometric functions and generalized basic hypergeometric functions $\;_{1}\Phi_{1}$, see also, \cite{Kumar2016,Kumar2017}. In \cite{sb} Srinivasulu and Bhagavan introduced the Irreducible representation of $SL(2, C)$ along with generating relations for generalized hypergeometric functions. Manocha \cite{ma} explored the Lie algebras associated with difference-differential operators and the Appell function $F_{1}$. Radulescu \cite{ra} studied of some special functions with Lie theory. In this paper we have introduced dynamical symmetry algebra of $\;_{1}F_{1}$ and $\Psi_{2}$ and by the induced group action arrived at certain identities for $\;_{1}F_{1}$ and $\Psi_{2}$ which in their turn lead to reduction formulae for hypergeometric functions of generating functions for different families $\;_{1}F_{1}$ and $\Psi_{2}$.

The confluent hypergeometric function $\;_{1}F_{1}(a;b;x)$ is defined by (see, \cite{sl1})

\begin{equation} \label{sec1eqn1}
\;_{1}F_{1}\left(a;b;x\right)= \sum_{s=0}^{\infty}\frac{\left(a\right)_{s}x^{s}}{s!\left(b\right)_{s}},
\end{equation}
where $a$, $b$, $x$ are complex numbers, $b$ is not a negative integer and $(a)_{s}$ is Pochhammer's symbol, defined as
\begin{eqnarray*}
\begin{split}
(a)_{s}=(a)(a+1)(a+2)...(a+s-1).
\end{split}
\end{eqnarray*}
The $\;_{1}F_{1}(a;b;x)$ satisfies the following differential recursion relations:
\begin{eqnarray}
\begin{split}
D\;_{1}F_{1}\left(a;b;x\right)=&\frac{a}{b}\;_{1}F_{1}\left(a+1;b+1;x\right);\;D=\frac{\mathrm{d}}{\mathrm{d}x},\\
\left(D-1\right)\;_{1}F_{1}\left(a;b;x\right) = &\frac{a-b}{b}\;_{1}F_{1}\left(a;b+1;x\right),\label{1.2}
\end{split}
\end{eqnarray}
\begin{eqnarray}
\begin{split}
\left(\Theta+a\right)\;_{1}F_{1}\left(a;b;x\right)=&a\;_{1}F_{1}\left(a+1;b;x\right);\Theta = x\frac{\mathrm{d}}{\mathrm{d}x},\\
\left(\Theta+b-1\right)\;_{1}F_{1}\left(a;b;x\right) = &\left(b-1\right)\;_{1}F_{1}\left(a;b-1;x\right),\label{1.3}
\end{split}
\end{eqnarray}
and
\begin{eqnarray}
\begin{split}
\left(\Theta+b-a-x\right)\;_{1}F_{1}\left(a;b;x\right) = &\left(b-a\right)\;_{1}F_{1}\left(a-1;b;x\right).\label{1.4}
\end{split}
\end{eqnarray}
\section{Dynamical Symmetry Algebra of $\;_{1}F_{1}$} \label{sec 2}
Let
\begin{eqnarray}
\begin{split}
f_{a,b}\left(x,y,z\right)=\frac{\Gamma\left(b-a\right)\Gamma(a)}{\Gamma(b)}\;_{1}F_{1}\left(a;b;x\right)y^{a}z^{b},\label{2.1}
\end{split}
\end{eqnarray}
be the basis elements of a subspace of analytical functions of variables $x$, $y$ and $z$ associated with basic confluent functions.

The $E$-operators for $\;_{1}F_{1}$ are
\begin{eqnarray}
\begin{split}
E_{a}=&y\left(\frac{\partial}{\partial\;x}+y\frac{\partial}{\partial\;y}\right),
\end{split}
\end{eqnarray} \label{2.2}
\begin{eqnarray}
\begin{split}
E_{a'} = &\frac{1}{y}\left(x\frac{\partial}{\partial\;x}-y\frac{\partial}{\partial\;y}+z\frac{\partial}{\partial\;z}-x\right),
\end{split}
\end{eqnarray}  \label{2.3}
\begin{eqnarray}
\begin{split}
E_{b}=&z\left(\frac{\partial}{\partial\;x}-1\right),
\end{split}
\end{eqnarray}  \label{2.4}
\begin{eqnarray}
\begin{split}
E_{b'}=&\frac{1}{z}\left(x\frac{\partial}{\partial\;x}+z\frac{\partial}{\partial\;z}-1\right),
\end{split}
\end{eqnarray} \label{2.5}
and
\begin{eqnarray}
\begin{split}
E_{ab}=&yz\frac{\partial}{\partial\;x}.
\end{split}
\end{eqnarray} \label{2.6}
The actions of theses $E$-operators on $f_{a,b}(x,y,z)$ is given by
\begin{eqnarray*}
\begin{split}
E_{a}f_{a,b}\left(x,y,z\right)=&y\left(x\frac{\partial}{\partial\;x}+y\frac{\partial}{\partial\;y}\right)f_{a,b}\left(x,y,z\right)\\
=&\left(x\frac{\partial}{\partial\;x}+a\right)\frac{\Gamma(b-a)\Gamma(a)}{\Gamma(b)}\;_{1}F_{1}\left(a;b;x\right)y^{a+1}z^{b}\\
=&\frac{\Gamma\left(b-a\right)\Gamma\left(a+1\right)}{\Gamma(b)}\;_{1}F_{1}\left(a+1;b;x\right)\,y^{a+1}z^{b}\\
=&\frac{\left(b-a-1\right)\Gamma\left(b-a-1\right)\Gamma\left(a+1\right)}{\Gamma(b)}\;_{1}F_{1}\left(a+1;b;x\right)\,y^{a+1}z^{b}\\
=&\left(b-a-1\right)f_{a+1,b}\left(x,y,z\right).
\end{split}
\end{eqnarray*}
Hence, we have
\begin{eqnarray}
\begin{split}
E_{a}f_{a,b}(x,y,z)=&(b-a-1)f_{a+1,b}(x,y,z).
\end{split}
\end{eqnarray} \label{2.7}
Similarly, we can derive the following results:
\begin{eqnarray}
\begin{split}
E_{a'}f_{a,b}(x,y,z)=&(a-1)f_{a-1,b}(x,y,z),
\end{split}
\end{eqnarray} \label{2.8}
\begin{eqnarray}
\begin{split}
E_{b}f_{a,b}(x,y,z)=&-f_{a,b+1}(x,y,z),
\end{split}
\end{eqnarray} \label{2.9}
\begin{eqnarray}
\begin{split}
E_{b'}f_{a,b}(x,y,z)=& (b-a-1)f_{a,b-1}(x,y,z),
\end{split}
\end{eqnarray} \label{2.10}
\begin{eqnarray}
\begin{split}
E_{ab}f_{a,b}(x,y,z)=f_{a+1,b+1}(x,y,z),
\end{split}
\end{eqnarray} \label{2.11}
\begin{eqnarray}
\begin{split}
I_{a}f_{a,b}(x,y,z)=af_{a,b}(x,y,z),
\end{split}
\end{eqnarray} \label{2.12}
\begin{eqnarray}
\begin{split}
I_{b}f_{a,b}(x,y,z)=bf_{a,b}(x,y,z),
\end{split}
\end{eqnarray} \label{2.13}
and
\begin{eqnarray}
\begin{split}
I f_{a,b}(x,y,z)=&f_{a,b}(x,y,z).
\end{split}
\end{eqnarray}   \label{2.14}
To find  $\exp(\alpha E_{a})$, we employ the operator
\begin{eqnarray*}
\begin{split}
E_{a}=&y\left(x\frac{\partial}{\partial\;x}+y\frac{\partial}{\partial\;y}\right),
\end{split}
\end{eqnarray*} \label{2.15}
with action
\begin{eqnarray*}
\begin{split}
E_{a}f_{a,b}(x,y,z)=&(b-a-1)f_{a+1,b}(x,y,z),
\end{split}
\end{eqnarray*}  \label{2.16}
for computing action of one parameter subgroup $\exp(\alpha E_{a})$ by usual multiplier representation theory we have to solve following differential equations:
\begin{eqnarray}
\begin{split}
\frac{\mathrm{d}y(\alpha)}{\mathrm{d}\alpha}=&y^{2}(\alpha),\;y(0)=y,
\end{split}
\end{eqnarray}  \label{2.17}
\begin{eqnarray*}
\begin{split}
\alpha=0,\, y(0)=y\Rightarrow \mathbf{k}=-\frac{1}{y},
\end{split}
\end{eqnarray*}
then
\begin{eqnarray}
\begin{split}
y(\alpha)=\frac{y}{1-\alpha y};
\end{split}
\end{eqnarray} \label{2.18}
and
\begin{eqnarray}
\begin{split}
\frac{\mathrm{d}x(\alpha)}{\mathrm{d}\alpha}=&x(\alpha)y(\alpha),\;x(0)=x
\end{split}
\end{eqnarray} \label{2.19}
and
\begin{eqnarray*}
\begin{split}
\alpha=&0,\, x(0)=x\Rightarrow \mathbf{k}=x,
\end{split}
\end{eqnarray*}
then
\begin{eqnarray}
\begin{split}
x(\alpha)=\frac{x}{1-\alpha y}.
\end{split}
\end{eqnarray} \label{2.20}
So that we have
\begin{eqnarray}
\begin{split}
\exp(\alpha\;E_{a})F_{a,b}(x,y,z)=\frac{\Gamma(b-a)\Gamma(a)}{\Gamma(b)}\;_{1}F_{1}\bigg{(}a;b;\frac{x}{1-\alpha y}\bigg{)}\bigg{(}\frac{y}{1-\alpha y}\bigg{)}^{a}z^{b}.
\end{split}
\end{eqnarray} \label{2.21}
On the other hand by direct expansion
\begin{eqnarray}
\begin{split}
&\exp(\alpha\;E_{a})F_{a,b}(x,y,z)=\sum_{\ell=0}^{\infty}\frac{\alpha^{\ell}}{\ell\;!}\bigg{(}E_{a}\bigg{)}^{\ell}F_{a,b}(x,y,z)\\
&=\sum_{\ell=0}^{\infty}\frac{\alpha^{\ell}\Gamma(b-a)}{\ell\;!\Gamma(b-a-\ell)}\frac{\Gamma(b-a-\ell)\Gamma(a+\ell)}{\Gamma(b)}\;_{1}F_{1}(a+\ell;b;x)y^{a+\ell}z^{b}.
\end{split}
\end{eqnarray}   \label{2.22}
Hence, we have
\begin{eqnarray}
\begin{split}
\left(1-\alpha y\right)^{-a}\;_{1}F_{1}\left(a;b;\frac{x}{1-\alpha y}\right)&=\sum_{\ell=0}^{\infty}\frac{(a)_{\ell}}{\ell\;!}\;_{1}F_{1}(a+\ell;b;x)\left(\alpha\;y\right)^{\ell}\;\;\;\; \left(\left|\alpha y\right|<1\right).
\end{split}
\end{eqnarray}  \label{2.23}
Setting $\alpha\;y\rightarrow \chi$, then we obtain
\begin{eqnarray}
\begin{split}
\left(1-\chi\right)^{-a}\;_{1}F_{1}\left(a;b;\frac{x}{1-\chi}\right)&=\sum_{\ell=0}^{\infty}\frac{(a)_{\ell}}{\ell!}\;_{1}F_{1}(a+\ell;b;x)\chi^{\ell}\;\;\;\; \left(\left|\chi\right|<1\right).
\end{split}
\end{eqnarray} \label{2.24}
Next, we use the operator
\begin{eqnarray*}
\begin{split}
E_{b}=&z\left(\frac{\partial}{\partial\;x}-1\right),
\end{split}
\end{eqnarray*}  \label{2.25}
with action
\begin{eqnarray*}
\begin{split}
E_{b}f_{a,b}(x,y,z)=&-f_{a,b+1}(x,y,z).
\end{split}
\end{eqnarray*} \label{2.26}

To compute the action of a one-parameter subgroup using the standard multiplier representation theory, we need to solve the corresponding differential equations.

To $\exp(\alpha E_{b})$
\begin{eqnarray}
\begin{split}
\frac{\mathrm{d}z(\alpha)}{\mathrm{d}\alpha}=&1,\,z(0)=z
\end{split}
\end{eqnarray} \label{2.27}
when
\begin{eqnarray*}
\begin{split}
\alpha=&0,\, z(0)=z\Rightarrow \mathbf{k}=z
\end{split}
\end{eqnarray*}
then
\begin{eqnarray}
\begin{split}
z(\alpha)=z+\alpha,
\end{split}
\end{eqnarray} \label{2.28}
and
\begin{eqnarray}
\begin{split}
\frac{\mathrm{d}x(\alpha)}{\mathrm{d}\alpha}=&\frac{x(\alpha)}{z(\alpha)},\, x(0)=x
\end{split}
\end{eqnarray} \label{2.29}
when
\begin{eqnarray*}
\begin{split}
\alpha=&0,\, x(0)=x\Rightarrow \mathbf{k}=\frac{x}{z}
\end{split}
\end{eqnarray*}
then
\begin{eqnarray}
\begin{split}
x(\alpha)=\frac{x}{z}(z+\alpha)
\end{split}
\end{eqnarray} \label{2.30}
and
\begin{eqnarray}
\begin{split}
\frac{\mathrm{d}u(\alpha)}{\mathrm{d}a}=&-\frac{z(\alpha)}{(\alpha)},
\end{split}
\end{eqnarray} \label{2.31}
\begin{eqnarray*}
\begin{split}
\alpha=&0,\,u(0)=1\Rightarrow \mathbf{k}=z
\end{split}
\end{eqnarray*}
then
\begin{eqnarray}
\begin{split}
u(\alpha)=\frac{z}{z+\alpha}.
\end{split}
\end{eqnarray} \label{2.32}
So that we have
\begin{eqnarray}
\begin{split}
&\exp(\alpha\;E_{b})F_{a,b}(x,y,z)\\
&=\frac{\Gamma(b-a)\Gamma(a)}{\Gamma(b)}z\;_{1}F_{1}(a+\ell;b+\ell;
\frac{x}{z}(z+\alpha))\left(z+\alpha\right)^{b-1}y^{a}.
\end{split}
\end{eqnarray} \label{2.33}
On the other hand by direct expansions, it yields
\begin{eqnarray}
\begin{split}
&\exp(\alpha\;E_{b})F_{a,b}(x,y,z)=\sum_{\ell=0}^{\infty}\frac{\alpha^{\ell}}{\ell!}\left(E_{b'}\right)^{\ell}F_{a,b}(x,y,z)\\
&=\sum_{\ell=0}^{\infty}\frac{(-1)^{\ell}}{\ell!}\;_{1}F_{1}(a;b+\ell;x)\left(\alpha z\right)^{\ell}y^{a}z^{b}.
\end{split}
\end{eqnarray} \label{2.34}
Equating two values of $\exp(\alpha\;E_{b})F_{a,b}(x,y,z)$, we arrive at the identity
\begin{eqnarray}
\begin{split}
&z\;_{1}F_{1}\left(a;b;\frac{x}{z}(z+\alpha)\right)y^{a}(z+\alpha)^{b-1}\\
&=\sum_{\ell=0}^{\infty}\frac{(b-a)_{\ell}}{\ell!(b)_{\ell}}\;_{1}F_{1}(a;b+\ell;x)\left(-\alpha z\right)^{\ell}y^{a}z^{b},
\end{split}
\end{eqnarray} \label{2.35}
which gives the generating relation
\begin{eqnarray}
\begin{split}
&\;_{1}F_{1}\left(a;b;\frac{x}{u}(z+\alpha)\right)\left(1+\frac{\alpha}{z}\right)^{b-1}\\
&=\sum_{\ell=0}^{\infty}\frac{(b-a)_{\ell}}{\ell!(b)_{\ell}}\;_{1}F_{1}(a;b+\ell;x)\left(-\alpha z\right)^{\ell}\;\;\;\left(\left|\frac{\alpha}{z}\right|<1\right).
\end{split}
\end{eqnarray} \label{2.36}
To find  $\exp(\alpha E_{a'})$, we employ the operator
\begin{eqnarray*}
\begin{split}
E_{a'}=&\frac{1}{y}\left(x\frac{\partial}{\partial\;x}-y\frac{\partial}{\partial\;y}+z\frac{\partial}{\partial\;z}-x\right),
\end{split}
\end{eqnarray*} \label{2.37}
with action
\begin{eqnarray*}
\begin{split}
E_{a'}f_{a,b}(x,y,z)=&(a-1)f_{a-1,b}(x,y,z),
\end{split}
\end{eqnarray*} 

\begin{eqnarray}
\begin{split}
\frac{\mathrm{d}y(\alpha)}{\mathrm{d}\alpha}=&-1,\,y(0)=y
\end{split}
\end{eqnarray} \label{2.38}
when
\begin{eqnarray*}
\begin{split}
\alpha=0,\,y(0)=y\Rightarrow \mathbf{k}=y,
\end{split}
\end{eqnarray*}
\begin{eqnarray}
\begin{split}
y(\alpha)=y-\alpha,
\end{split}
\end{eqnarray}  \label{2.39}
\begin{eqnarray}
\begin{split}
\frac{\mathrm{d}x(\alpha)}{\mathrm{d}\alpha}=&\frac{x(\alpha)(1-x(\alpha))}{y(\alpha)},
\end{split}
\end{eqnarray} \label{2.40}
when
\begin{eqnarray*}
\begin{split}
\alpha=0,\,x(0)=x\Rightarrow \mathbf{k}=\frac{xy}{1-x},
\end{split}
\end{eqnarray*}
then
\begin{eqnarray}
\begin{split}
x(\alpha)=\frac{xy}{y-\alpha(1-x)},
\end{split}
\end{eqnarray} \label{2.41}
\begin{eqnarray}
\begin{split}
\frac{\mathrm{d}z(\alpha)}{\mathrm{d}\alpha}=&-\frac{z(\alpha)x(\alpha)}{y(\alpha)},
\end{split}
\end{eqnarray} \label{2.42}
\begin{eqnarray*}
\begin{split}
\alpha=0,\, z(0)=z\Rightarrow \mathbf{k}=z,
\end{split}
\end{eqnarray*}
\begin{eqnarray}
\begin{split}
z(\alpha)=\frac{z(y-\alpha)}{y-\alpha(1-x)},
\end{split}
\end{eqnarray} \label{2.43}
\begin{eqnarray}
\begin{split}
&\exp(\alpha\;E_{a'})f_{a,b}(x,y,z)\\
&=\frac{\Gamma(b-a)\Gamma(a)}{\Gamma(b)}\;_{1}F_{1}\left(a;b;\frac{xy}{y-\alpha(1-x)}\right)(y-\alpha)^{a}\left(\frac{z(y-\alpha)}{(y-\alpha(1-x))}\right)^{b},
\end{split}
\end{eqnarray} \label{2.44}
\begin{eqnarray}
\begin{split}
&\exp(\alpha\;E_{a'})f_{a,b}(x,y,z)=\sum_{\ell=0}^{\infty}\frac{\alpha^{\ell}}{\ell!}\left(E_{a'}\right)^{\ell}\,f_{a,b}(x,y,z)\\
&=\sum_{\ell=0}^{\infty}\frac{\alpha^{\ell}(a-\ell)_{\ell}}{\ell!}\frac{\Gamma(b-a+\ell)\Gamma(a-\ell)}{\Gamma(b)}\;_{1}F_{1}(a-\ell;b;x)y^{a-\ell}z^{b}.
\end{split}
\end{eqnarray} \label{2.45}
Hence, we have
\begin{eqnarray}
\begin{split}
&\;_{1}F_{1}\left(a;b;\frac{xy}{y-\alpha(1-x)}\right)\left(\frac{(y-\alpha)}{(y-\alpha(1-x))}\right)^{b}\left(1-\frac{\alpha}{y}\right)^{a}\\
&=\sum_{\ell=0}^{\infty}\frac{(b-a)_{\ell}}{\ell!}\;_{1}F_{1}(a-\ell;b;x)\left(\frac{\alpha}{y}\right)^{\ell}\;\;\;\; \left(\left|\frac{\alpha}{y}\right|<1,\, \left|\frac{\alpha(1-x)}{y}\right|<1\right).
\end{split}
\end{eqnarray} \label{2.46}
To find  $\exp(\alpha E_{b'})$, we utilize the following operator:
\begin{eqnarray*}
\begin{split}
E_{b'}=&\frac{1}{z}\left(x\frac{\partial}{\partial\;x}+z\frac{\partial}{\partial\;z}-1\right),
\end{split}
\end{eqnarray*} \label{2.47}
with action
\begin{eqnarray*}
\begin{split}
E_{b'}f_{a,b}(x,y,z)=& (b-a-1)f_{a,b-1}(x,y,z),
\end{split}
\end{eqnarray*} 

\begin{eqnarray}
\begin{split}
\frac{\mathrm{d}z(\alpha)}{\mathrm{d}\alpha}=&1,
\end{split}
\end{eqnarray} \label{2.48}
when
\begin{eqnarray*}
\begin{split}
\alpha=0,\,z(0)=z\Rightarrow \mathbf{k}=z,
\end{split}
\end{eqnarray*}
\begin{eqnarray}
\begin{split}
z(\alpha)=z+\alpha,
\end{split}
\end{eqnarray}  \label{2.49}
\begin{eqnarray}
\begin{split}
\frac{\mathrm{d}x(\alpha)}{\mathrm{d}\alpha}=&\frac{x(\alpha)}{z(\alpha)},\,x(0)=x,
\end{split}
\end{eqnarray}  \label{2.50}
\begin{eqnarray*}
\begin{split}
\alpha=0,\,x(0)=x\Rightarrow \mathbf{k}=\frac{x}{z},
\end{split}
\end{eqnarray*}
then
\begin{eqnarray}
\begin{split}
x(\alpha)=\frac{x(z+\alpha)}{z}.
\end{split}
\end{eqnarray} \label{2.51}
So that, we have
\begin{eqnarray}
\begin{split}
\exp(\alpha\;E_{b'})f_{a,b}(x,y,z)=\frac{\Gamma(b-a)\Gamma(a)}{\Gamma(b)}\;_{1}F_{1}\left(a;b;\frac{x(z+\alpha)}{z}\right)y^{a}(z+\alpha)^{b}.
\end{split}
\end{eqnarray} \label{2.52}
Conversely, using directed expansion, we obtain
\begin{eqnarray}
\begin{split}
&\exp(\alpha\;E_{b'})f_{a,b}(x,y,z)=\sum_{\ell=0}^{\infty}\frac{\alpha^{\ell}}{\ell!}\left(E_{b'}\right)^{\ell}f_{a,b}(x,y,z)\\
&=\sum_{\ell=0}^{\infty}\frac{\alpha^{\ell}(b-a-\ell)_{\ell}}{\ell!}\frac{\Gamma(b-a-\ell)\Gamma(a)}{\Gamma(b-\ell)}\;_{1}F_{1}(a;b-\ell;x)y^{a}z^{b}.
\end{split}
\end{eqnarray} \label{2.53}
which gives the generating relation
\begin{eqnarray}
\begin{split}
&\;_{1}F_{1}\left(a;b;x(1+\frac{\alpha}{z})\right)\left(1+\frac{\alpha}{z}\right)^{b}\\
&=\sum_{\ell=0}^{\infty}\frac{(b-\ell)_{\ell}}{\ell !}\;_{1}F_{1}(a;b-\ell,d;x)\left(\frac{\alpha}{z}\right)^{\ell}\;\; \left(\left|\frac{\alpha}{z}\right|<1\right).
\end{split}
\end{eqnarray} \label{2.54}
Setting $\frac{\alpha}{z}\rightarrow \chi$, then we obtain
\begin{eqnarray}
\begin{split}
&\;_{1}F_{1}\left(a;b;x(1+\chi)\right)\left(1+\chi\right)^{b}\\
&=\sum_{\ell=0}^{\infty}\frac{(b-\ell)_{\ell}}{\ell !}\;_{1}F_{1}(a;b-\ell,d;x)\chi^{\ell}\;\; \left(\left|\chi\right|<1\right).
\end{split}
\end{eqnarray} \label{2.55}

Now, we use the operator
\begin{eqnarray}
\begin{split}
E_{ab}=& yz\frac{\partial}{\partial\;x},
\end{split}
\end{eqnarray} \label{2.56}
with the action given by
\begin{eqnarray}
\begin{split}
E_{ab}f_{a,b}(x,y,z)=&f_{a+1,b+1}(x,y,z).
\end{split}
\end{eqnarray} \label{2.57}
Computing action of one parameter subgroup by usual multiplier representation theory we solve these differential equations
\begin{eqnarray}
\begin{split}
\frac{\mathrm{d}x(\alpha)}{\mathrm{d}\alpha}=&yz,\,x(0)=x
\end{split}
\end{eqnarray} \label{2.58}
when
\begin{eqnarray*}
\begin{split}
\alpha=&0,\,x(0)=x\Rightarrow \mathbf{k}=x
\end{split}
\end{eqnarray*}
then
\begin{eqnarray}
\begin{split}
x(\alpha)=x+\alpha zy.
\end{split}
\end{eqnarray} \label{2.59}
Thus, we have
\begin{eqnarray}
\begin{split}
\exp(\alpha\;E_{ab})F_{a,b}(x,y,z)=\frac{\Gamma(b-a)\Gamma(a)}{\Gamma(b)}\;_{1}F_{1}(a;b;x+\alpha yz)y^{a}z^{b}.
\end{split}
\end{eqnarray} \label{2.60}
On the other hand, through direct expansion, we obtain
\begin{eqnarray}
\begin{split}
&\exp(\alpha\;E_{ab})F_{a,b}(x,y,z)=\sum_{\ell=0}^{\infty}\frac{\alpha^{\ell}}{\ell !}\left(E_{ab}\right)^{\ell}F_{a,b}(x,y,z)\\
&=\sum_{\ell=0}^{\infty}\frac{1}{\ell !}\frac{\Gamma(b-a)\Gamma(a+\ell)}{\Gamma(b+\ell)}\;_{1}F_{1}(a+\ell;b+\ell;x)\left(\alpha\;y\;z\right)^{\ell}y^{a}z^{b}.
\end{split}
\end{eqnarray} \label{2.61}
By equating the two values of $\exp(\alpha\;E_{ab})F_{a,b}(x,y,z)$, we get
\begin{eqnarray}
\begin{split}
\;_{1}F_{1}(a;b;x+u\alpha z)y^{a}z^{b}=\sum_{\ell=0}^{\infty}\frac{(a)_{\ell}}{\ell !(b)_{\ell}}\;_{1}F_{1}(a+\ell;b+\ell;x)\left(\alpha yz\right)^{\ell}y^{a}z^{b},
\end{split}
\end{eqnarray} \label{2.62}
which gives the generating relation, given as
\begin{eqnarray}
\begin{split}
\;_{1}F_{1}(a;b;x+\alpha yz)=\sum_{\ell=0}^{\infty}\frac{(a)_{\ell}}{\ell !(b)_{\ell}}\;_{1}F_{1}(a+\ell;b+\ell;x)\left(\alpha yz\right)^{\ell}.
\end{split}
\end{eqnarray} \label{2.63}
By setting $\alpha\;yz\rightarrow \chi$, we obtain
\begin{eqnarray}
\begin{split}
\;_{1}F_{1}(a;b;x+\chi)=\sum_{\ell=0}^{\infty}\frac{(a)_{\ell}}{\ell !(b)_{\ell}}\;_{1}F_{1}(a+\ell;b+\ell;x)\chi^{\ell}.
\end{split}
\end{eqnarray} \label{2.64}

\section{The Dynamical Symmetry Algebra of $\Psi_{2}(a;b,c;x)$} \label{sec 3}

The Humbert function $\Psi_{2}(a;b,c;x)$ is defined by (see \cite{emot})
\begin{eqnarray}
\begin{split}
\Psi_{2}(a;b,c;x)=&\sum_{m,n=0}^{\infty}\frac{(a)_{m+n}x^{m}y^{n}}{m!n!(b)_{m}(c)_{n}}=\sum_{m,n=0}^{\infty}\frac{(a)_{m}(a+m)_{n}x^{m}y^{n}}{m!n!(b)_{m}(c)_{n}}\\
=&\sum_{m=0}^{\infty}\frac{(a)_{m}x^{m}}{m!(b)_{m}}\;_{1}F_{1}(a+m;c;y),
\end{split}
\end{eqnarray} \label{3.1}
or we can express it as
\begin{eqnarray}
\begin{split}
\Psi_{2}(a;b,c;x,y)=\sum_{n=0}^{\infty}\frac{(a)_{n}y^{n}}{n!(c)_{n}}\;_{1}F_{1}(a+n;b;x).
\end{split}
\end{eqnarray}   \label{3.2}
Let be set the new Humbert functions
\begin{eqnarray}
\begin{split}
\Psi_{a,b,c}(x,y,z,u,t)=\Psi_{2}(a;b,c;x,y)z^{a}u^{b}t^{c},
\end{split}
\end{eqnarray} \label{3.3}
be the basic elements of a subspace of analytical functions for three variables $x,y,z$ associated with hypergeometric function $\;_{1}F_{1}$ defined as $\Psi_{2}$.

The dynamical symmetry algebra of $\Psi_{2}(a;b,c;x,y)$ is a complex Lie algebra generated by $E$-operators, which are referred to as raising or lowering operators due to their ability to raise or lower the corresponding suffix

\subsection{E-operators}
The $E$-operators for $\Psi_{2}(a;b,c;x,y)$ are
\begin{eqnarray}
\begin{split}
E_{a}=&z\left(x\frac{\partial}{\partial\;x}+y\frac{\partial}{\partial\;y}+z\frac{\partial}{\partial\;z}\right),
\end{split}
\end{eqnarray} \label{3.4}

\begin{eqnarray}
\begin{split}
E_{b}=&\frac{1}{u}\left(x\frac{\partial}{\partial\;x}+u\frac{\partial}{\partial\;u}-1\right),
\end{split}
\end{eqnarray} \label{3.5}
\begin{eqnarray}
\begin{split}
E_{c}=&\frac{1}{t}\left(y\frac{\partial}{\partial\;y}+t\frac{\partial}{\partial\;t}-1\right),
\end{split}
\end{eqnarray} \label{3.6}
\begin{eqnarray}
\begin{split}
E_{a,b}=&uz\frac{\partial}{\partial\;x},
\end{split}
\end{eqnarray} \label{3.7}
and
\begin{eqnarray}
\begin{split}
E_{a,c}=&zt\frac{\partial}{\partial\;y}.
\end{split}
\end{eqnarray} \label{3.8}
The action of these $E$-operators on $\Psi_{a,b,c}(x,y,z,u,t)$ is given by
\begin{eqnarray*}
\begin{split}
&E_{a}\Psi_{a,b,c}(x,y,z,u,t)=E_{a}\sum_{n=0}^{\infty}\frac{(a)_{n}y^{n}}{n!(c)_{n}}\;_{1}F_{1}(a+n;b;x)z^{a}u^{b}t^{c}\\
=&z\left(x\frac{\partial}{\partial\;x}\sum_{n=0}^{\infty}\frac{(a)_{n}y^{n}}{n!(c)_{n}}\;_{1}F_{1}(a+n;b;x)+y\frac{\partial}{\partial\;y}\sum_{n=0}^{\infty}\frac{(a)_{n}y^{n}}{n!(c)_{n}}\;_{1}F_{1}(a+n;b;x)\right.\\
&\left.+z\frac{\partial}{\partial\;z}\sum_{n=0}^{\infty}\frac{(a)_{n}y^{n}}{n!(c)_{n}}\;_{1}F_{1}(a+n;b;x)\right)z^{a}u^{b}t^{c}\\
=&z\left(\sum_{n=0}^{\infty}\frac{(a)_{n}y^{n}}{n!(c)_{n}}\left((a+n)\;_{1}F_{1}(a+n+1;b;x)-(a+n)\;_{1}F_{1}(a+n;b;x)\right)z^{a}u^{b}t^{c}\right.\\
&\left. +\sum_{n=0}^{\infty}\frac{n(a)_{n}y^{n}}{n!(c)_{n}}\;_{1}F_{1}(a+n;b;x)z^{a}u^{b}t^{c}+a\sum_{n=0}^{\infty}\frac{(a)_{n}y^{n}}{n!(c)_{n}}\;_{1}F_{1}(a+n;b;x)z^{a}u^{b}t^{c}\right)\\
=&z\left(\sum_{n=0}^{\infty}\frac{(a)_{n}y^{n}}{n!(c)_{n}}(a+n)\;_{1}F_{1}(a+n+1;b;x)-\sum_{n=0}^{\infty}\frac{(a)_{n}y^{n}}{n!(c)_{n}}(a+n)\;_{1}F_{1}(a+n;b;x)\right.\\
&\left. +\sum_{n=0}^{\infty}\frac{(a+n)(a)_{n}y^{n-1}}{n!(c)_{n}}\;_{1}F_{1}(a+n;b;x)\right)z^{a}u^{b}t^{c}\\
&=z\sum_{n=0}^{\infty}\frac{(a+n)(a)_{n}y^{n}}{n!(c)_{n}}\;_{1}F_{1}(a+n+1;b;x)z^{a}u^{b}t^{c}\\
=&\sum_{n=0}^{\infty}\frac{a(a+1)_{n}y^{n}}{n!(c)_{n}}\;_{1}F_{1}(a+n+1;b;x)z^{a+1}u^{b}t^{c}=a\Psi_{a+1,b,c}(x,y,z,u,t).
\end{split}
\end{eqnarray*}
Hence, we get
\begin{eqnarray}
\begin{split}
&E_{a}\Psi_{a,b,c}(x,y,z,u,t)=a\Psi_{a+1,b,c}(x,y,z,u,t).
\end{split}
\end{eqnarray} \label{3.9}
By using the above techniques, we obtain
\begin{eqnarray}
\begin{split}
E_{b}\Psi_{a,b,c}(x,y,z,u,t)=&(b-1)\Psi_{a,b-1,c}(x,y,z,u,t),
\end{split}
\end{eqnarray} \label{3.10}
\begin{eqnarray}
\begin{split}
E_{c}\Psi_{a,b,c}(x,y,z,u,t)=&(c-1)\Psi_{a,b,c-1}(x,y,z,u,t),
\end{split}
\end{eqnarray} \label{3.11}
\begin{eqnarray}
\begin{split}
E_{a,b}\Psi_{a,b,c,d}(x,y,s,u,v,t)=\frac{a}{b}\Psi_{a+1,b+1,c}(x,y,z,u,t),
\end{split}
\end{eqnarray} \label{3.12}
and
\begin{eqnarray}
\begin{split}
E_{a,c}\Psi_{a,b,c,d}(x,y,s,u,v,t)=\frac{a}{c}\Psi_{a+1;b,c+1}(x,y,s,u,v,t).
\end{split}
\end{eqnarray} \label{3.13}

The upper factor in each bracket should be linked with a plus sign, while the lower factor corresponds to a minus sign. These $E$-operators, along with the three maintenance operators $I_{a}$, $I_{b}$, and $I_{c}$ as well as the identity operators $I$, collectively form a basis. We have
\begin{eqnarray}
\begin{split}
I_{a}=&z\frac{\partial}{\partial{z}},\\
I_{b}=&u\frac{\partial}{\partial{u}},\\
I_{c}=&t\frac{\partial}{\partial{t}},\\
I=&1.
\end{split}
\end{eqnarray} \label{3.14}
\begin{eqnarray}
\begin{split}
I_{a}\Psi_{a,b,c}(x,y,z,u,t)=&a\Psi_{a,b,c}(x,y,z,u,t),\\
I_{b}\Psi_{a,b,c}(x,y,z,u,t)=&b\Psi_{a,b,c}(x,y,z,u,t),\\
I_{c}\Psi_{a,b,c}(x,y,z,u,t)=&c\Psi_{a,b,c}(x,y,z,u,t).
\end{split}
\end{eqnarray}  \label{3.15}
This section presents a group-theoretic framework for deriving reduction formulas for the Humbert function in five variables.

To find  $\exp(\alpha E_{a})$, we employ the operator
\begin{eqnarray*}
\begin{split}
E_{a}=z\left(x\frac{\partial}{\partial x}+y\frac{\partial}{\partial y}+z\frac{\partial}{\partial z}\right),
\end{split}
\end{eqnarray*}
with action
\begin{eqnarray*}
\begin{split}
E_{a}\Psi_{a,b,c}(x,y,z,u,t)=&a\Psi_{a+1,b,c}(x,y,z,u,t).
\end{split}
\end{eqnarray*}

To compute the action of the one-parameter subgroup $\exp(\alpha E_{a})$ using the standard multiplier representation theory, we have to solve following differential equations:
\begin{eqnarray}
\begin{split}
\frac{dz(\alpha)}{d\alpha}=&z^{2}(\alpha),\; z(0)=z
\end{split}
\end{eqnarray}  \label{3.16}
\begin{eqnarray*}
\begin{split}
\alpha=0,\;z(0)=z\Rightarrow \mathbf{k}=-\frac{1}{z},
\end{split}
\end{eqnarray*}
then
\begin{eqnarray}
\begin{split}
z(\alpha)=\frac{z}{1-\alpha z},
\end{split}
\end{eqnarray}  \label{3.17}
and
\begin{eqnarray}
\begin{split}
\frac{dx(\alpha)}{d\alpha}=&x(\alpha)z(\alpha),\;x(0)=x,
\end{split}
\end{eqnarray}  \label{3.18}
then
\begin{eqnarray*}
\begin{split}
\alpha=&0,\, x(0)=x\Rightarrow \mathbf{k}=x,
\end{split}
\end{eqnarray*}
\begin{eqnarray}
\begin{split}
x(\alpha)=\frac{x}{1-\alpha z},
\end{split}
\end{eqnarray} \label{3.19}
\begin{eqnarray}
\begin{split}
\frac{\mathrm{d}y(\alpha)}{\mathrm{d}\alpha}=&y(\alpha)z(\alpha),\;y(0)=y,
\end{split}
\end{eqnarray} \label{3.20}
\begin{eqnarray*}
\begin{split}
\alpha=&0,\, y(0)=y\Rightarrow \mathbf{k}=y,
\end{split}
\end{eqnarray*}
\begin{eqnarray}
\begin{split}
y(\alpha)=\frac{y}{1-\alpha z},
\end{split}
\end{eqnarray} \label{3.21}
So that, we get
\begin{eqnarray}
\begin{split}
\exp(\alpha\;E_{a})\Psi_{a,b,c}(x,y,z,u,t)=\Psi_{2}\left(a;b,c;\frac{x}{1-\alpha z},\frac{y}{1-\alpha z}\right)\left(\frac{z}{1-\alpha z}\right)^{a}u^{b}t^{c}.
\end{split}
\end{eqnarray} \label{3.22}
On the other hand by direct expansion, we have
\begin{eqnarray}
\begin{split}
&\exp(\alpha\;E_{a})\Psi_{a,b,c}(x,y,z,u,t)=\sum_{\ell=0}^{\infty}\frac{\alpha^{\ell}}{\ell !}\left(E_{a}\right)^{\ell}\Psi_{a,b,c}(x,y,z,u,t)\\
&=\sum_{\ell=0}^{\infty}\frac{\alpha^{\ell}(a)_{\ell}}{\ell !(b)_{\ell}}\Psi_{a+\ell,b,c}(x,y,z,u,t)\\
&=\sum_{\ell=0}^{\infty}\frac{\alpha^{\ell}(a)_{\ell}}{\ell !(b)_{\ell}}\Psi_{2}(a+\ell;b,c;x,y,z,u,t)z^{a+\ell}u^{b}t^{c}\\
&=\sum_{\ell=0}^{\infty}\frac{(a)_{\ell}}{\ell !(b)_{\ell}}\Psi_{2}(a+\ell;b,c;x,y,z,u,t)\left(\alpha\;z\right)^{\ell}z^{a}u^{b}t^{c}.
\end{split}
\end{eqnarray} \label{3.23}
Hence, we arrive at
\begin{eqnarray}
\begin{split}
&\left(1-\alpha z\right)^{-a}\Psi_{2}\left(a;b,c;\frac{x}{1-\alpha z},\frac{y}{1-\alpha z}\right)\\
&=\sum_{\ell=0}^{\infty}\frac{(a)_{\ell}}{\ell!(b)_{\ell}}\Psi_{2}(a+\ell;b,c;x,y,z,u,t)\left(\alpha\;z\right)^{\ell}.
\end{split}
\end{eqnarray} \label{3.24}
By comparing the two values, we establish the following:
\begin{eqnarray}
\begin{split}
&(1-\alpha z)^{-a}\Psi_{2}\left(a;b,c;\frac{x}{1-\alpha z},\frac{y}{1-\alpha z}\right)z^{a}u^{b}t^{c}=\sum_{\ell=0}^{\infty}\frac{\alpha^{\ell}}{\ell!}\left(E_{a}\right)^{\ell}\Psi_{a,b,c}(x,y,z,u,t)\\
&=\sum_{\ell=0}^{\infty}\frac{\alpha^{\ell}(a)_{\ell}}{\ell!}\Psi_{a+\ell,b+\ell,c}(x,y,z,u,t)\\
&=\sum_{\ell=0}^{\infty}\frac{\alpha^{\ell}(a)_{\ell}}{\ell!}\sum_{m,n=0}^{\infty}\frac{(a+\ell)_{m+n}x^{m}y^{n}}{m!n!(b)_{m}(c)_{n}}z^{a+\ell}u^{b}t^{c}\\
&=\sum_{\ell,m,n=0}^{\infty}\frac{\alpha^{\ell}(a)_{\ell}}{\ell!}\frac{(a+\ell)_{m+n}x^{m}y^{n}}{m!n!(b)_{m}(c)_{n}}z^{a+\ell}u^{b}t^{c},
\end{split}
\end{eqnarray} \label{3.25}
which after simplification gives
\begin{eqnarray}
\begin{split}
&(1-\alpha z)^{-a}\Psi_{2}\left(a;b,c;\frac{x}{1-\alpha z},\frac{y}{1-\alpha z}\right)\\
&=\sum_{\ell,m,n=0}^{\infty}\frac{(a)_{\ell+m+n}x^{m}y^{n}}{\ell!m!n!(b)_{m}(c)_{n}}(\alpha z)^{\ell}\;\;\; \left(\left|\alpha z\right|<1\right).
\end{split}
\end{eqnarray} \label{3.26}
By setting $\alpha\;z\rightarrow \chi$ and in view of the following definition:
\begin{eqnarray}
\begin{split}
\Psi_{2}\left(a;b,c;x,y,z\right)=\sum_{\ell,m,n=0}^{\infty}\frac{(a)_{\ell+m+n}}{(b)_{m}(c)_{n}}\frac{x^{m}y^{n}z^{\ell}}{\ell!m!n!},
\end{split}
\end{eqnarray} \label{3.27}
we arrive at the reduction formula, given by
\begin{eqnarray}
\begin{split}
\Psi_{2}\left(a;b,c;x,y,\chi\right)=(1-\chi)^{-a}\Psi_{2}\left(a;b,c;\frac{x}{1-\chi},\frac{y}{1-\chi}\right)\;\;\; \left(\left|\chi\right|<1\right).
\end{split}
\end{eqnarray} \label{3.28}
Next, we use the operator
\begin{eqnarray*}
\begin{split}
E_{b}=&\frac{1}{u}\left(x\frac{\partial}{\partial\;x}+u\frac{\partial}{\partial\;u}-1\right),
\end{split}
\end{eqnarray*}
with action
\begin{eqnarray*}
\begin{split}
E_{b}\Psi_{a,b,c}(x,y,z,u,t)=&(b-1)\Psi_{a,b-1,c}(x,y,z,u,t).
\end{split}
\end{eqnarray*}
To compute the action of the one-parameter subgroup using the standard multiplier representation theory, we have to solve the following differential equations:
To $\exp(\alpha E_{b})$
\begin{eqnarray}
\begin{split}
\frac{\mathrm{d}u(\alpha)}{\mathrm{d}\alpha}=&1,\;u(0)=u,
\end{split}
\end{eqnarray} \label{3.29}
when
\begin{eqnarray*}
\begin{split}
\alpha=&0,u(0)=u\Rightarrow \mathbf{k}=u,
\end{split}
\end{eqnarray*}
then
\begin{eqnarray}
\begin{split}
u(\alpha)=u+\alpha,
\end{split}
\end{eqnarray} \label{3.30}
and
\begin{eqnarray}
\begin{split}
\frac{\mathrm{d}x(\alpha)}{\mathrm{d}\alpha}=&\frac{x(\alpha)}{u(\alpha)},\;x(0)=x,
\end{split}
\end{eqnarray} \label{3.31}
when
\begin{eqnarray*}
\begin{split}
\alpha=&0,\, x(0)=x\Rightarrow \mathbf{k}=\frac{x}{u},
\end{split}
\end{eqnarray*}
then
\begin{eqnarray}
\begin{split}
x(\alpha)=\frac{x}{u}(u+\alpha),
\end{split}
\end{eqnarray} \label{3.32}
and
\begin{eqnarray}
\begin{split}
\frac{\mathrm{d}v(\alpha)}{\mathrm{d}a}=&-\frac{v(\alpha)}{u(\alpha)},\;z(0)=z,
\end{split}
\end{eqnarray} \label{3.33}
when
\begin{eqnarray*}
\begin{split}
\alpha=&0,v(0)=1\Rightarrow \mathbf{k}=u,
\end{split}
\end{eqnarray*}
then
\begin{eqnarray}
\begin{split}
v(\alpha)=\frac{u}{u+\alpha}.
\end{split}
\end{eqnarray} \label{3.34}
So that we obtain
\begin{eqnarray}
\begin{split}
\exp(\alpha E_{b'})\Psi_{a,b,c}(x,y,z,u,t)=u\Psi_{2}(a+\ell;b+\ell,c;\frac{x}{u}(u+\alpha),y)\left(u+\alpha\right)^{b-1}z^{a}t^{c}.
\end{split}
\end{eqnarray} \label{3.35}
On the other hand by direct expansions, it yields
\begin{eqnarray}
\begin{split}
&\exp(\alpha\;E_{b})\Psi_{a,b,c}(x,y,z,u,t)=\sum_{\ell=0}^{\infty}\frac{\alpha^{\ell}}{\ell!}\left(E_{b'}\right)^{\ell}\Psi_{a,b,c}(x,y,z,u,t)\\
&=\sum_{\ell=0}^{\infty}\frac{(b-\ell)_{\ell}}{\ell!}\Psi_{2}(a;b-\ell,c;x,y)\left(\frac{\alpha}{u}\right)^{\ell}z^{a}u^{b}t^{c}.
\end{split}
\end{eqnarray} \label{3.36}
By equating two values of $\exp(\alpha\;E_{b})\Psi_{a,b,c}(x,y,z,u,t)$, we obtain the following identity:
\begin{eqnarray}
\begin{split}
&u\Psi_{2}\left(a;b,c;\frac{x}{u}(u+\alpha),y\right)z^{a}(u+\alpha)^{b-1}t^{c}\\
&=\sum_{\ell=0}^{\infty}\frac{(b-\ell)_{\ell}}{\ell!}\Psi_{2}(a;b-\ell,c;x,y)\left(\frac{\alpha}{u}\right)^{\ell}z^{a}u^{b}t^{c},
\end{split}
\end{eqnarray} \label{3.37}
which gives the generating relation
\begin{eqnarray}
\begin{split}
&\Psi_{2}\left(a;b,c;\frac{x}{u}(u+\alpha),y\right)\left(1+\frac{\alpha}{u}\right)^{b-1}\\
&=\sum_{\ell=0}^{\infty}\frac{(b-\ell)_{\ell}}{\ell!}\Psi_{2}(a;b-\ell,c;x,y)\left(\frac{\alpha}{u}\right)^{\ell}\;\;\; \left(\left|\frac{\alpha}{u}\right|<1\right).
\end{split}
\end{eqnarray} \label{3.38}
Next, we use the operator
\begin{eqnarray*}
\begin{split}
E_{c}=&\frac{1}{t}\left(y\frac{\partial}{\partial\;y}+t\frac{\partial}{\partial\;t}-1\right),
\end{split}
\end{eqnarray*}
with action
\begin{eqnarray*}
\begin{split}
E_{c}\Psi_{a,b,c}(x,y,z,u,t)=&(c-1)\Psi_{a,b,c-1}(x,y,z,u,t).
\end{split}
\end{eqnarray*}
Computing action of one parameter subgroup by usual multiplier representation theory, we have to solve following differential equations:
\begin{eqnarray}
\begin{split}
\frac{\mathrm{d}t(\alpha)}{\mathrm{d}\alpha}=&1,
\end{split}
\end{eqnarray} \label{3.39}
when
\begin{eqnarray*}
\begin{split}
\alpha=&0,t(0)=t\Rightarrow \mathbf{k}=t,
\end{split}
\end{eqnarray*}
then
\begin{eqnarray}
\begin{split}
t(\alpha)=t+\alpha,
\end{split}
\end{eqnarray} \label{3.40}
\begin{eqnarray}
\begin{split}
\frac{\mathrm{d}y(\alpha)}{\mathrm{d}\alpha}=&\frac{y(\alpha)}{t(\alpha)},\;x(0)=x,
\end{split}
\end{eqnarray} \label{3.41}
when
\begin{eqnarray*}
\begin{split}
\alpha=&0,\, y(0)=y\Rightarrow \mathbf{k}=\frac{y}{t},
\end{split}
\end{eqnarray*}
then
\begin{eqnarray}
\begin{split}
y(\alpha)=\frac{y}{t}(t+\alpha),
\end{split}
\end{eqnarray} \label{3.42}
\begin{eqnarray}
\begin{split}
\frac{\mathrm{d}v(\alpha)}{\mathrm{d}a}=&-\frac{v(\alpha)}{t(\alpha)},\, z(0)=z,
\end{split}
\end{eqnarray} \label{3.43}
when
\begin{eqnarray*}
\begin{split}
\alpha=&0,\, v(0)=1\Rightarrow \mathbf{k}=t
\end{split}
\end{eqnarray*}
then
\begin{eqnarray}
\begin{split}
v(\alpha)=\frac{t}{t+\alpha}.
\end{split}
\end{eqnarray} \label{3.44}
So that we have
\begin{eqnarray}
\begin{split}
\exp(\alpha\;E_{c})\Psi_{a,b,c}(x,y,z,u,t)=t\Psi_{2}(a+\ell;b+\ell,c;x,\frac{y}{t}(t+\alpha))\left(t+\alpha\right)^{c-1}z^{a}u^{b}.
\end{split}
\end{eqnarray} \label{3.45}
On the other hand by direct expansions, it gives
\begin{eqnarray}
\begin{split}
&\exp(\alpha\;E_{c'})\Psi_{a,b,c}(x,y,z,u,t)=\sum_{\ell=0}^{\infty}\frac{\alpha^{\ell}}{\ell!}\left(E_{c'}\right)^{\ell}\Psi_{a,b,c}(x,y,z,u,t)\\
&=\sum_{\ell=0}^{\infty}\frac{(c-\ell)_{\ell}}{\ell!}\Psi_{2}(a;b,c-\ell;x,y)\left(\frac{\alpha}{t}\right)^{\ell}z^{a}u^{b}t^{c}.
\end{split}
\end{eqnarray} \label{3.46}
By equating two values of $\exp(\alpha\;E_{c})\Psi_{a,b,c}(x,y,z,u,t)$, we arrive at the identity
\begin{eqnarray}
\begin{split}
&t\Psi_{2}\left(a;b,c;x,\frac{y}{t}(t+\alpha)\right)z^{a}(t+\alpha)^{c-1}u^{b}\\
&=\sum_{\ell=0}^{\infty}\frac{(c-\ell)_{\ell}}{\ell!}\Psi_{2}(a;b,c-\ell;x,y)\left(\frac{\alpha}{t}\right)^{\ell}z^{a}u^{b}t^{c},
\end{split}
\end{eqnarray} \label{3.47}
which gives the generating relation
\begin{eqnarray}
\begin{split}
&\Psi_{2}\left(a;b,c;x,\frac{y}{t}(t+\alpha)\right)\left(1+\frac{\alpha}{t}\right)^{c-1}\\
&=\sum_{\ell=0}^{\infty}\frac{(c-\ell)_{\ell}}{\ell!}\Psi_{2}(a;b,c-\ell;x,y)\left(\frac{\alpha}{t}\right)^{\ell}\;\;\; \left(\left|\frac{\alpha}{t}\right|<1\right).
\end{split}
\end{eqnarray} \label{3.48}
Now, we use the operator
\begin{eqnarray*}
\begin{split}
E_{a,b}=&zu\frac{\partial}{\partial\;x},
\end{split}
\end{eqnarray*}
with action
\begin{eqnarray*}
\begin{split}
E_{a,b}\Psi_{a,b,c}(x,y,z,u,t)=&\frac{a}{b}\Psi_{a+1,b+1,c}(x,y,z,u,t).
\end{split}
\end{eqnarray*}
To compute the action of the one-parameter subgroup using the standard multiplier representation theory, we solve the following differential equations:
\begin{eqnarray}
\begin{split}
\frac{\mathrm{d}x(\alpha)}{\mathrm{d}\alpha}=&zu,\;x(0)=x,
\end{split}
\end{eqnarray} \label{3.49}
when
\begin{eqnarray*}
\begin{split}
\alpha=&0,x(0)=x\Rightarrow \mathbf{k}=x,
\end{split}
\end{eqnarray*}
then
\begin{eqnarray}
\begin{split}
x(\alpha)=x+u\alpha z.
\end{split}
\end{eqnarray} \label{3.50}
Thus
\begin{eqnarray}
\begin{split}
\exp(\alpha E_{ab})\Psi_{a,b,c}(x,y,z,u,t)=\Psi_{2}(a;b,c;x+u\alpha z,y)z^{a}u^{b}t^{c}.
\end{split}
\end{eqnarray} \label{3.51}
Alternatively, through direct expansion, we obtain:
\begin{eqnarray}
\begin{split}
&\exp(\alpha E_{ab})\Psi_{a,b,c}(x,y,z,u,t)=\sum_{\ell=0}^{\infty}\frac{\alpha^{\ell}}{\ell !}\left(E_{ab}\right)^{\ell}\Psi_{a,b,c}(x,y,z,u,t)\\
&=\sum_{\ell=0}^{\infty}\frac{(a)_{\ell}}{\ell!(b)_{\ell}}\Psi_{2}(a+\ell;b+\ell,c;x,y)\left(\alpha\;z\;u\right)^{\ell}z^{a}u^{b}t^{c}.
\end{split}
\end{eqnarray} \label{3.52}
Equating the two values of $\exp(\alpha\;E_{ab})\Psi_{a,b,c}(x,y,z,u,t)$, we get
\begin{eqnarray}
\begin{split}
\Psi_{2}(a;b,c;x+u\alpha z,y)z^{a}u^{b}t^{c}=\sum_{\ell=0}^{\infty}\frac{(a)_{\ell}}{\ell\;!(b)_{\ell}}\Psi_{2}(a+\ell;b+\ell,c;x,y)\big{(}zu\alpha\big{)}^{\ell}z^{a}u^{b}t^{c},
\end{split}
\end{eqnarray} \label{3.53}
which gives the generating relation
\begin{eqnarray}
\begin{split}
\Psi_{2}(a;b,c;x+u\alpha z,y)=\sum_{\ell=0}^{\infty}\frac{(a)_{\ell}}{\ell!(b)_{\ell}}\Psi_{2}(a+\ell;b+\ell,c;x,y)\left(zu\alpha\right)^{\ell}.
\end{split}
\end{eqnarray} \label{3.54}
Setting $u\alpha z \rightarrow \chi$, we get the generating relation
\begin{eqnarray}
\begin{split}
\Psi_{2}(a;b,c;x+\chi,y)=\sum_{\ell=0}^{\infty}\frac{(a)_{\ell}}{\ell!(b)_{\ell}}\Psi_{2}(a+\ell;b+\ell,c;x,y)\chi^{\ell}.
\end{split}
\end{eqnarray} \label{3.55}
finally, we use the operator
\begin{eqnarray*}
\begin{split}
E_{a,c}=&zt\frac{\partial}{\partial\;y}
\end{split}
\end{eqnarray*}
with action
\begin{eqnarray*}
\begin{split}
E_{a,c}\Psi_{a,b,c}(x,y,z,u,t)=&\frac{a}{c}\Psi_{a+1,b,c+1}(x,y,z,u,t).
\end{split}
\end{eqnarray*}
Computing action of one parameter subgroup by usual multiplier representation theory, we solve following differential equations:
\begin{eqnarray}
\begin{split}
\frac{\mathrm{d}y(\alpha)}{\mathrm{d}\alpha}=&zt,\;y(0)=y,
\end{split}
\end{eqnarray} \label{3.56}
when
\begin{eqnarray*}
\begin{split}
\alpha=&0,y(0)=y\Rightarrow \mathbf{k}=y,
\end{split}
\end{eqnarray*}
then
\begin{eqnarray}
\begin{split}
y(\alpha)=y+zt\alpha.
\end{split}
\end{eqnarray} \label{3.57}
Thus, we have
\begin{eqnarray}
\begin{split}
\exp(\alpha\;E_{ac})\Psi_{a,b,c}(x,y,z,u,t)=\Psi_{2}(a;b,c;x,y+zt\alpha)z^{a}u^{b}t^{c}.
\end{split}
\end{eqnarray}\label{3.58}
On the other hand by direct expansion, we get
\begin{eqnarray}
\begin{split}
&\exp(\alpha\;E_{ac})\Psi_{a,b,c}(x,y,z,u,t)=\sum_{\ell=0}^{\infty}\frac{\alpha^{\ell}}{\ell!}\left(E_{ac}\right)^{\ell}\Psi_{a,b,c}(x,y,z,u,t)\\
&=\sum_{\ell=0}^{\infty}\frac{(a)_{\ell}}{\ell!(c)_{\ell}}\Psi_{2}(a+\ell;b,c+\ell;x,y)\left(\alpha\;z\;t\right)^{\ell}z^{a}u^{b}t^{c}.
\end{split}
\end{eqnarray} \label{3.59}
By equating the two values of $\exp(\alpha\;E_{ac})\Psi_{a,b,c}(x,y,z,u,t)$, we get
\begin{eqnarray}
\begin{split}
\Psi_{2}(a;b,c;x,y+zt\alpha)z^{a}u^{b}t^{c}=\sum_{\ell=0}^{\infty}\frac{(a)_{\ell}}{\ell!(c)_{\ell}}\Psi_{2}(a+\ell;b,c+\ell;x,y)\big{(}zt\alpha\big{)}^{\ell}z^{a}u^{b}t^{c},
\end{split}
\end{eqnarray} \label{3.60}
which gives the generating relation
\begin{eqnarray}
\begin{split}
\Psi_{2}(a;b,c;x,y+zt\alpha)=\sum_{\ell=0}^{\infty}\frac{(a)_{\ell}}{\ell!(c)_{\ell}}\Psi_{2}(a+\ell;b,c+\ell;x,y)\left(zt\alpha\right)^{\ell}.
\end{split}
\end{eqnarray} \label{3.61}
By setting $zt\alpha \rightarrow \chi$, we get the generating relation
\begin{eqnarray}
\begin{split}
\Psi_{2}(a;b,c;x,y+\chi)=\sum_{\ell=0}^{\infty}\frac{(a)_{\ell}}{\ell !(c)_{\ell}}\Psi_{2}(a+\ell;b,c+\ell;x,y)\chi^{\ell}.
\end{split}
\end{eqnarray} \label{3.62}
The symmetry algebras corresponding to the different confluent forms of the $\;_{1}F_{1}(a;b;x)$ and $\Psi_{2}$ Humbert functions can likewise be calculated using the aforementioned methods. As an alternative, these Lie algebras can be found as contractions of the previously calculated symmetry algebras.
\section{Conclusion and concluding remarks} \label{sec 4}
The dynamical symmetry algebras of the confluent hypergeometric function  $\;_{1}F_{1}(a;b;x)$ and the $\Psi_{2}$ Humbert function can be easily constructed using our method, as it requires no complex calculations. A key result of our analysis is the explicit proof that different generalized hypergeometric functions transform according to representations of specific Lie algebras. The findings presented in this study are intended to have broad applications across various fields, including mathematics, physics, statistics, and engineering.

\noindent \textbf{Declarations}\\
The authors have no conflicts of interest to declare that are relevant to the content of this article.

\noindent \textbf{Availability of data and materials}\\
No data were used for this study.

\noindent \textbf{Competing interests}\\
The authors declare that they have no competing interests.

\noindent \textbf{Funding}\\
No funding was received for conducting this study.

\end{document}